\documentclass{article}

\usepackage[utf8]{inputenc}

\usepackage[left=1in,right=1in,top=1in,bottom=1in,includefoot,includehead,headheight=13.6pt]{geometry}

\usepackage{xcolor}
\definecolor{darkolivegreen}{rgb}{0.33, 0.42, 0.18}
\definecolor{celestialblue}{rgb}{0.29, 0.59, 0.82}

\usepackage[colorlinks,
            linkcolor=darkolivegreen,
            citecolor=darkolivegreen,
            urlcolor=darkolivegreen]{hyperref}

\hypersetup{pdftitle={CGTJS}, pdfauthor={Felipe Galarce Marin}, 
}

\usepackage{fullpage}
\usepackage{amsmath, amsthm, amssymb}
\usepackage{thmtools, thm-restate}
\usepackage{dsfont}
\usepackage{mathtools}
\mathtoolsset{showonlyrefs=true}
\usepackage{graphicx}   
\usepackage{subfigure}
\usepackage{color}
\usepackage{setspace}
\theoremstyle{definition}

\numberwithin{equation}{section} % Number eqs with section-number

\newcommand{\bR}{\ensuremath{\mathbb{R}}}

\def\[{\left[}
\def\]{\right]}
\def\<{\langle}
\def\>{\rangle}
\def\({\left(}
\def\){\right)}
\def\[{\left [}
\def\]{\right]}
\def\({\left(}
\def\){\right)}

\newcommand{\norm}[1]{\Vert #1 \Vert}

\usepackage{enumitem}

\usepackage{acronym,authblk}
\usepackage{todonotes}
\usepackage{booktabs}
\usepackage{mathtools} % coloneqq
\usepackage{natbib}
\usepackage{isomath}
\renewcommand{\vec}{\vectorsym}
\newcommand{\mat}{\matrixsym}

\newcommand{\R}{\mathbb{R}}
\newcommand{\bu}{\boldsymbol{u}}

\newcommand{\revNB}[1]{\textcolor{black}{#1}}

\acrodef{dof}[DoF]{Degree of Freedom}
\acrodefplural{dof}[DoFs]{Degrees of Freedom}
\acrodef{AMG}{Algebraic Multigrid}
\acrodef{CG}{Conjugate Gradient}
\acrodef{FE}{Finite Element}
\acrodef{AA}{Anderson Acceleration}
\acrodef{LS}{Least-Squares}
\acrodef{DE}{Differential Equation}
\acrodefplural{DE}[DEs]{Differential Equations}
\acrodef{ODE}{Ordinary Differential Equation}
\acrodefplural{ODE}[ODEs]{Ordinary Differential Equations}
\acrodef{PDE}{Partial Differential Equation}
\acrodefplural{PDE}[PDEs]{Partial Differential Equations}

\begin{document}

\title{Accelerated implicitization: Robust fixed-point iterations arising from an explicit scheme}

\author[1]{Nicolas A. Barnafi}
\author[2]{Felipe Galarce}
\author[3]{Pablo Brubeck}

\affil[1]{Instituto de Ingeniería Matemática y Computacional. Pontificia Universidad Católica de Chile, Chile.}
\affil[2]{School of Civil Engineering. Pontificia Universidad Católica de Valparaíso, Chile}
\affil[3]{Mathematical Institute. University of Oxford, Oxford UK.}
\date{\today}

\maketitle

\begin{abstract}

\noindent \textbf{Purpose:} In this work we propose a general strategy for solving the possibly non-linear problems arising from implicit time discretizations as a sequence of explicit solutions. This sequence can be as unstable as the base explicit discretization, which can be improved through Anderson acceleration.

\noindent \textbf{Design/methodology/approach:} We propose using explicit fixed-point sub-iterations for non linear problems combined with the Anderson Acceleration (AA) technique to improve convergence and speed, verifying its usability and scalability in three non-linear differential equations.

We provide an error analysis to establish the expected properties of the proposed strategy, both for time and space-time based problems. Through several examples, we show how simple it is to setup this method and then expose its strengths by doing parameter sensitivity.

\noindent \textbf{Findings:} The proposed method is simple to implement and yields acceptable performance in a wide range of problems. This type of method is best suited for situations in which matrix assembly can be prohibitively expensive, or whenever a good preconditioner for the implicit problem is out of reach, such as in highly convective fluid flows. \revNB{In tests where non-accelerated iterations are convergent, acceleration provides an overall reduction in the iterations of up to a 96\%.}

\noindent \textbf{Originality:} This work formalizes the delay of implicit terms in an implicit time discretization, provides a succinct error analysis and further enhances the strategy through Anderson acceleration. The results are encouraging and well founded in existing theory, thus establishing the foundations for interesting future research.

%[TODO]
%In this work we propose using explicit time discretizations to obtain arbitrary time discretizations upon convergence through fixed-point subiterations. An error analysis shows that the convergence is related to the Lipschitz constant of the function defining the dynamics, and that a CFL-type condition arises for 

%This work propose a general pipeline to solve fixed-point problems that typically arise when discretizing partial differential equations in time. We rely on explicit Picard iterations enhanced with the already well-known Anderson acceleration method. We therefore iterate on the fixed-point loop at a negligible cost (involving sometimes only the inverse of a mass matrix), while enlarging the convergence ratio with the acceleration scheme, allowing larger time-steps and both, lower computational complexity (\textcolor{red}{to be cheked}) and lower wall-clock CPU times compared to fully implicit or semi-implicit fixed-point iterations.
\end{abstract}


%\tableofcontents

\section{Introduction}
Time discretization methods are a broad topic that ranges over all applications that have time dependence, which are naturally abundant \citep{temam2012infinite}. Within the community interested in using these methods to numerically approximate equations, many topics are of interest, such as the development of high-order approximation schemes, stability, and structure preservation. See \citep{hairer1993solving} for an introductory presentation on all of the aforementioned topics, among other interesting ones.

One lingering topic that is transversal to all of these questions is whether the resulting scheme is \emph{explicit} or \emph{implicit}, i.e. if a linear system needs to be solved at each time instant. Most structure preserving schemes are actually implicit \citep{wanner1996solving}, and there is an intuition in the community that acknowledges that \emph{implicitness} is the price to pay for structure in general, despite there being some exceptions in Hamiltonian systems. Still, it is very uncommon to see explicit schemes being used in elastodynamics or fluid mechanics, for example, as they suffer from numerical instabilities. Motivated by this, in this work we built upon the following idea: Given an implicit time discretization, we use as a solution strategy the fixed point iteration induced by lagging the evaluation of the implicitly-treated terms. The resulting sub-problem is equivalent to the problem obtained when using an explicit discretization. An error analysis shows that this is typically as unstable as the original explicit discretization, but an acceleration strategy renders the overall method much more robust. We refer to the process of leveraging explicit iterations to solve implicit ones as \emph{implicitization}. We highlight that upon convergence, the solution obtained matches that of the implicit discretization by construction.

The acceleration strategy we adopt is \ac{AA} \citep{anderson1965iterative}, which has gained a lot of traction due to its solid theoretical guarantees. As a small review: it typically accelerates linearly convergent sequences up to quadratic convergence \citep{toth2015convergence,evans2020proof} with an increased radius of convergence \citep{lupo2019convergence}. If one accelerates a Richardson iteration for linear systems, the resulting accelerated sequence matches that of GMRES in exact arithmetic \citep{walker2011anderson}, and a similar procedure for nonlinear systems yields a multi-secant method \citep{walker2011anderson}. \ac{AA} computes an optimal combination of residuals, and thus its main computational overhead is that of storing such residuals and then solving a \ac{LS} problem. A memory width, known as depth in the context of \ac{AA}, plays the role of restart in GMRES, with the difference of it being a moving window that reduces the storage required. The cost of the \ac{LS} solution can be reduced by performing acceleration only every $k$ fixed-point iterations, known as \emph{alternation} \citep{suryanarayana2019alternating}, and also by doing approximate solutions through sketching \citep{lupo2024anderson}, which can be tailored to specific physics as well \citep{barnafi2025two}. 

\revNB{\ac{AA} has been widely used in different applications and problem structures. In direct problems, it has been very successfully used to accelerate Picard iterations in Navier-Stokes both in Newtonian \citep{pollock2019anderson} and non-Newtonian regimes \citep{diaz2024accelerated,pollock2023anderson}, as well as similar schemes for multiphysics problems such as Boussinesq \citep{pollock2021acceleration} and Magnetohydrodynamics \citep{dong2024finite}. It has also been used to solve Electrophysiology equations \citep{barnafi2024robust}, and in inverse problems for proximal point methods \citep{mai2020anderson} and image registration \citep{barnafi2025tree}. Other methodological works include accelerating descent methods for Reinforcement and Deep Learning models \citep{geist2018anderson,pasini2022anderson}.}

\revNB{The scope of this work is the proposal of a simple and robust iterative framework that leverages \ac{AA} to obtain the solution of an implicit (possibly nonlinear) discretization scheme through sub-iterations of an explicit (linear) one. To maintain the scope of the work contained, we will focus on Euler discretization schemes and show analytical constraints for the unaccelerated subiterations that are greatly alleviated by \ac{AA}.} To do this, we first formulate the implicitization method rigorously and provide a convergence analysis in Section~\ref{sec:method}. We will establish that this strategy can be interpreted as a Quasi-Newton method. In Section~\ref{sec:anderson} we show how the method can be improved using Anderson acceleration, and verify how this method provides only a small computational overhead, mainly acceptable for \acp{PDE}, and show the complete resulting algorithm. We conclude with a number of numerical tests to validate our methodology in \acp{ODE} and \acp{PDE}.

\section{Implicitization formulation}\label{sec:method}

To define implicitization, we will consider an \ac{ODE} in $\R$ given by finding a function $x:[0,T]\to \R^d$ that solves a Cauchy problem
    \begin{equation}\label{eq:ode}
    \left\{
        \begin{aligned}
            \dot x(t) &= f(t, x),\\
            x(0)   &= x_0,
        \end{aligned}
    \right.
    \end{equation}
with $f:(0,T)\times \R^d \to \R^d$, $T>0$ and for some given initial condition $x_0$ in $\R^d$. To discretize \eqref{eq:ode} one can consider a set of equidistant points $0=t^1<t^2<\dots < t^N = T$ such that their increments define a timestep $\Delta t \coloneqq t^{n+1} - t^n > 0$, and then define approximate collocation points $x^n\approx x(t^n)$ which become the new variables.  An implicit approximation of \eqref{eq:ode} yields the nonlinear system
    \begin{equation}\label{eq:ode-imp}
        x^n - x^{n-1} = \Delta t f(t^n, x^n),
    \end{equation}
For this scheme, we consider the implicitization procedure as the fixed-point iteration given by the sequence $(x^{n,\ell})_\ell$, defined by delaying the nonlinear term in \eqref{eq:ode-imp} as 
    \begin{equation}\label{eq:ode-imp-iter}
        x^{n,\ell} - x^{n-1} = \Delta t f(t^n, x^{n,\ell-1}).
    \end{equation}
Naturally, one should consider $x^{n,0}=x^{n-1}$, and note that if $x^{n,\ell}$ converges as a sequence in $\ell$, it converges to the solution of the implicit formulation \eqref{eq:ode-imp}, i.e. $x^{n,\ell} \to x^n$. In that regard, we note that an implicitization of an explicit scheme is simply a solution strategy for the (possibly) nonlinear problem that \eqref{eq:ode-imp} gives rise to. Indeed, this strategy can be equivalently rewritten as a quasi-Newton method. To note this, we rewrite \eqref{eq:ode-imp-iter} as an equation for the increment $\delta x^{n,\ell} \coloneqq x^{n,\ell} - x^{n, \ell - 1}$, which yields the system
    $$ \delta x^{n,\ell} = -(x^{n,\ell} - x^{n-1}) + \Delta t f(t^n, x^{n,\ell}). $$
The right hand side is the residual expression of \eqref{eq:ode-imp}, given by  $ F^n(x^\star) \coloneqq (x^\star - x^{n-1}) - \Delta t f(t_n, x^\star)$. Everything together yields the iteration system 
    \begin{equation}\label{eq:qn}
        \delta x^{n,\ell} = -F^n(x^{n,\ell})
        \end{equation}
for each $\ell$. This establishes the claimed equivalence, for an approximation matrix given simply by the identity. One can readily see that in a time dependent \ac{PDE}, the matrix approximating the Jacobian would be a mass matrix, which is inexpensive to solve, and can be robustly preconditioned 
by its diagonal \citep{wathen1987realistic}, or further approximated
through a lumping process \citep{hinton1976note}. 

\subsection{Convergence analysis}
For simplicity, we will look at the case in which the function $f$ is Lipschitz continuous on its second argument, i.e. there is a constant $L>0$ such that  for all $t\in [0,T]$ it holds that
    $$|f(t, x^1) - f(t, x^2)| \leq L|x^1 - x^2|.$$ 
In that case, consider the difference of two consecutive iterates of \eqref{eq:ode-imp-iter}, which using the Lipschitz continuity of $f$ gives
    $$ |x^{n,\ell+1} - x^{n,\ell}| = \Delta t|f(t^n,x^{n,\ell}) - f(t^n, x^{n, \ell -1})| \leq \Delta t L|x^{n,\ell} - x^{n, \ell-1}|.$$
We thus obtain that the implicitization sequence converges for a sufficiently small time step, i.e.
    $$ \Delta t \leq \frac 1 L.$$
Bounds relating the timestep to the inverse of the Lipschitz constant are common in stability analysis \citep{hairer1993solving}, and show that a large Lipschitz constant directly relates to the difficulty of the proposed methodology for a specific problem.

\subsection{On the spatial dependency of implicitization}
\label{sec:example-heat}

In this section we see how to apply this procedure to problems with spatial dependency, and to fix ideas we look at the heat equation:
    \begin{equation}\label{eq:heat}
        \left\{
        \begin{aligned}
            \dot u - \Delta u &= 0 && \text{in $(0,T)\times \Omega$} \\
            u &= 0 && \text{on $(0,T)\times \partial \Omega$} \\
            u &= u_0 && \text{on $\{0\}\times \Omega$},
        \end{aligned}
        \right.
    \end{equation}
whose implicit discretization in time with a \ac{FE} discretization is given by: Given an initial field $u_h^0$, find $\{u_h^n\}_n$ in a discrete space $V_h$ (which we consider to be $H^1$-conforming piecewise linear polynomials) such that
    \begin{equation}\label{eq:heat-imp}
        (u_h^{n} - u_h^{n-1}, v_h)  - \Delta t(\nabla u_h^n, \nabla v_h) = 0 \qquad\forall v_h \in V_h.
    \end{equation}
The implicitization iteration is given by the sequence $\{u_h^{n,\ell}\}_\ell$ induced by the problem
    \begin{equation}\label{eq:heat-imp-iter}
        (u_h^{n,\ell} - u_h^{n-1}, v_h)  - \Delta t(\nabla u_h^{n,\ell-1}, \nabla v_h) = 0 \qquad\forall v_h \in V_h.
    \end{equation}
Given the error $e_h^{n,\ell} \coloneqq u_h^{n,\ell} - u_h^{n,\ell-1}$ we get the error equation
    \begin{equation}\label{eq:error}
        (e_h^{n,\ell}, v_h) + \Delta t (\nabla e_h^{n,\ell-1}, \nabla v_h) = 0 \qquad \forall v_h\in V_h.
    \end{equation}
Recall the inverse inequality \citep{ern2004theory} for some generic mesh-independent constant $C>0$
    $$ \| v_h \|_1 \leq \frac C h \| v_h \|_0\qquad \forall v_h \in V_h,$$
and set $v_h=e_h^{n,\ell}$ in \eqref{eq:error} to obtain through Young's inequality that
    $$\| e_h^{n,\ell}\|_0^2 \leq \frac {C\Delta t} 2 \|e_h^{n,\ell}\|_1^2 +  \frac{C\Delta t} 2 \| e_h^{n,\ell-1} \|_1^2 \Rightarrow  \left(1  - \frac{C\Delta t}{2 h^2}\right)\|e_h^{n,\ell}\|_0^2 \leq \frac {C\Delta t}{2 h^2} \| e_h^{n,\ell-1} \|_0^2.$$
    To ensure that $\| e_h^{n,\ell} \| \leq \|e_h^{n,\ell-1}\|$ we thus require the following CFL condition
    $$ \Delta t \leq \tilde C h^2$$
for some mesh-independent constant $C$. We conclude that the fixed-point iteration that arises from implicitization is guaranteed to converge only under the same conditions that the explicit time discretization requires for stability. We leave open the question whether this fixed-point iteration is tightly related to time instabilities, and focus on proposing a computationally inexpensive way to render this iterative process more robust. We present the tools to do this in Section~\ref{sec:anderson}.

\section{Fixed-point iterations and Anderson acceleration}\label{sec:anderson}
We first briefly introduce how \ac{AA} works. For this, consider a fixed-point iteration of the form  $ u_{k+1} = G(u_{k})$,  which naturally arises from the numerical solution of non-linear equations.  Let us denote $r_k \coloneqq u_k - G(u_k) \in \R^N$, the $k$-th residual. Additionally, let $u^1, \ldots u^m$ be the last $m_k \coloneqq \min \{ k, m \}$ iterations ($m$ is the \emph{depth}).  We look for an optimal combination of the previous vectors, given by the solution of the following problem:
\begin{equation}\label{eq:anderson}
    \begin{aligned}
        \min_{(\alpha_0,\ldots,\alpha_{m_k}) \in \bR^{{m_k}+1}} \left\| \sum_{i=0}^{m_k} \alpha_i r_i\right\|_2^2 \\
        \text{s.t. } \sum_{i=0}^{m_k} \alpha_i = 1.
    \end{aligned}
\end{equation}
This formulation, while very explicit, is not usually implemented in practice due to the presence of a constraint \citep{walker2011anderson}. We convert it into an unconstrained problem by setting $\alpha_0 = 1 - \sum_{i=1}^m \alpha_i$, which yields:

\begin{equation}\label{eq:ls}
\min_{(\alpha_0,\ldots,\alpha_{m_k}) \in \bR^{{m_k}+1}} \norm{ \sum_{i=0}^{m_k} \alpha_i r_i }  =  \min_{(\alpha_1,\ldots,\alpha_{m_k}) \in \bR^{m_k}} \norm{ \sum_{i=1}^{m_k} \alpha_i (r_i - r_0) + r_0 },
\end{equation}

which is a linear least-squares problem for $\vec{\alpha} = (\alpha_1, \ldots, \alpha_{m_k})$, whose normal equations are given by:
$$
\Delta^T \Delta \vec{\alpha} = \Delta^T r_0,
$$
with $\Delta \in \bR^{N \times m_k}$, a matrix whose $i$-th column is given by $r_i - r_0$. The problem can be solved at minimal cost using a QR factorization, avoiding assembling $\Delta^T \Delta$, therefore avoiding squaring the condition number of the normal equations.

In this case, Anderson acceleration can be interpreted as a multi-secant method \citep{walker2011anderson}. Then, the overall solution method can be seen as an extension of the previously defined quasi-Newton method, but where a multi-secant update is performed to enrich the initial Jacobian approximation. For a robust formulation, see \citep[Eq. (3.1)]{walker2011anderson}.

\revNB{At each iteration, the method requires the fixed point map $G$ to be evaluated anyway. Having this in consideration, \ac{AA} has a small but non-negligible computational overhead with respect to the evaluation of $G$. The memory overhead is the $\Delta$ matrix in $\R^{n\times m_k}$ which stores the past $m_k$ residual differences, and the performance overhead is given by the solution of the least-squares problem \eqref{eq:ls}. Using a QR factorization, it can be solved in linear complexity $O(nm_k)$, and thus the algorithm remains linear in the problem size. Efforts have been done to leverage Randomized Linear Algebra algorithms in order to subsample the rows of the least squares problem and make this step even less expensive \citep{lupo2024anderson,barnafi2025two}. Naturally, the case $m=0$ reduces to the Picard iteration given by $u_{k+1} = G(u_k)$.}

\revNB{
\subsection{Accelerated implicitized iterations}\label{sec:accel-imp}
In this section we show how implicitized sub-iterations can be accelerated with \ac{AA} in the context of \acp{PDE}. To this end, consider the differential equation $\dot{u} = \mathcal{L}u$ given in terms of a differential operator $\mathcal L: V\to V'$ formulated in a distributional sense. Let $F(u; v) = \langle \mathcal{L} u, v \rangle$. The weak formulation is to find $u$ in $V$ such that
    $$ (\dot{u}, v) = F(u; v)  \quad\forall\, v\in V, $$
with can be implicitly discretized with a timestep $\Delta t$ as 
    $$ (u^{n}, v) = (u^{n-1}, v) + \Delta t F(u^{n}; v)  \quad\forall\, v\in V. $$
    The implicitized discrete sub-iteration is given by considering a discrete space $V_h \subset V$ and an initial guess $u_h^{n,0} = u_h^{n-1}$ and then finding $u_h^{n,k}$ for each $k$ such that:
    $$ (u_h^{n,k}, v_h) = (u_h^{n-1}, v_h) + \Delta t F(u_h^{n,k-1}; v_h)  \quad \forall\, v_h \in V_h. $$
    Defining the basis functions $\{\varphi_i\}_{i=1}^N$ that span $V_h$, the mass matrix $\mat M_{ij} \coloneqq (\varphi_j, \varphi_i)$, the right hand side vector $\vec F_i(u_h^{n,k-1}) \coloneqq  F(u_h^{n,k-1}; \varphi_i )$, and the vector of coefficients $\vec U^{\eta, \xi}$ such that $u_h^{\eta,\xi} = \sum_{i=1}^N U^{\eta,\xi}_i \varphi_i$, we obtain the discrete problem to be solved for $\vec U^{n,k}$: 
    \begin{equation} 
            \mat M \vec U^{n,k} = \mat M  \vec U^{n-1} + \Delta t\vec F(\vec U^{n,k-1}).
    \end{equation}
Thus, in the notation of Section~\ref{sec:anderson}, we obtain the fixed point operator at time $n$ given by the solution map
    $$ G(\vec U) \coloneqq \vec U^{n-1} + \Delta t \mat M^{-1}\vec F(\vec U), $$
    which can then be used for \ac{AA}. For example, under the CFL condition obtained in Section~\ref{sec:example-heat} the operator $G$ induced by $F=\Delta$ is a linearly convergent contraction. The contractive property of the fixed-point operator implies that the rate of convergence to a fixed point is linear. In addition, it has been shown in \citep{evans2020proof} that when \ac{AA} is applied to a linearly converging sequence, it generates an accelerated sequence that converges no more than quadratically. Thus, thanks to the contractive property, the accelerated sequence is also convergent \citep{toth2015convergence}, and the linear convergence rate means that the fixed-point operator $G$ is amenable to acceleration.
}

\section{Numerical Examples}

In this section we present several numerical tests that validate our proposed methodology. For this we do the following:

    \begin{enumerate}
        \item A convergence test on an electrophysiology \ac{ODE}. 
        \item Implicitization for the heat equation
        \item Implicitization for the Navier-Stokes equations for a more challenging test.
    \end{enumerate}

    \begin{figure}[!htbp]
        \centering
        \includegraphics[width=0.7\linewidth]{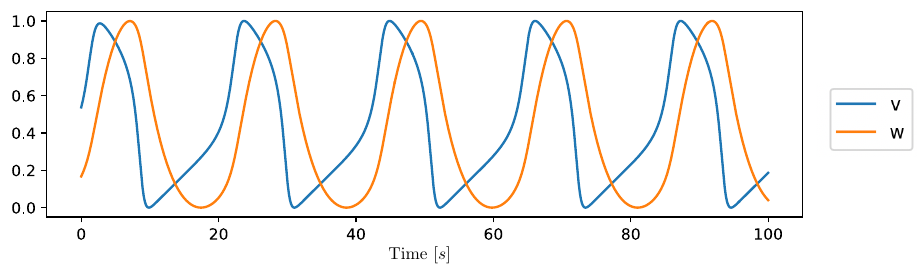}
        \caption{Normalized numerical solutions to \eqref{eq:fgn} using fully implicit scheme. The parameter configuration is set to $a=0.7$, $b=0.8$, $\tau=5$, $R=1$ and $I_\text{ext}=0.5$.}
        \label{fig:ground_truth_FNG}
    \end{figure}

\subsection{Electrophysiology model for cell-membrane potential}

The electric activity within connective tissue is usually modeled macroscopically using the mono-domain or bi-domain models \citep{ColliFranzone2014}, which aim to capture the propagation of the electric potential in muscle cells. These models are reaction-diffusion equations where the reaction terms account for the local cell membrane potential response. There are several models to account for the cell activation, action potential, and refractory behavior \citep{BENSON202160}. To test out our numerical approach, we select one of the simplest yet physiologically consistent models: the FitzHugh–Nagumo equations. The cell dynamics are governed by two main unknowns: the cell electric potential $v : [0,T] \rightarrow \bR$ ($T>0$), and an internal recovery variable $w : [0,T] \rightarrow \bR$. 

    \begin{figure}[!htbp]
        \centering
        \includegraphics[width=0.8\linewidth]{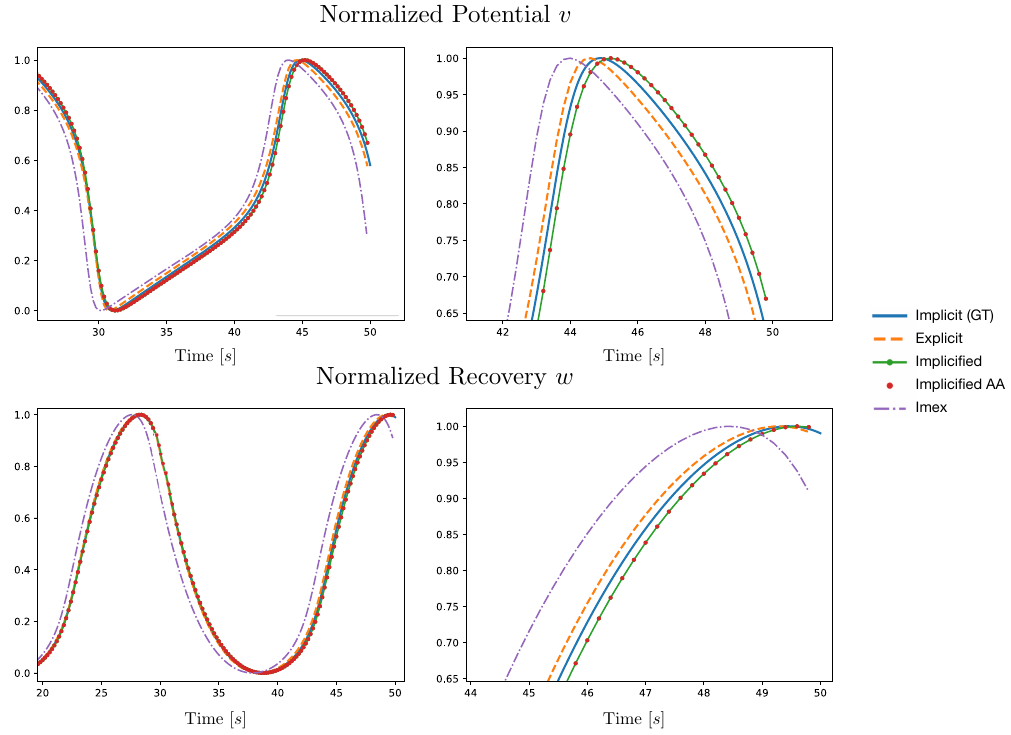}
        \caption{Numerical solutions with $\Delta t = 0.1$.}
        \label{fig:normalized_vw}
    \end{figure}

The system of ordinary differential equations reads:
\begin{equation}
\left\{\begin{aligned}
\dot{v} &= v - \frac{v^3}{3} - w + R I_{\text{ext}}, & \quad (0,T) \\
\tau \dot{w} &= v + a - b w, & \quad (0,T),
\end{aligned}\right.
\label{eq:fgn}
\end{equation}
closed with some given initial conditions. The parameters $a$, $b$, $R$, $\tau$ and $I_\text{ext}$ are real constants.

For our comparisons, we will consider the following discretization alternatives:
\begin{enumerate}
    \item Fully implicit approach (ground truth): Here we set a Picard approach with implicit terms avoiding only the cubic non-linearity in the first conservation equation from \eqref{eq:fgn}. This approach boils down to solve: 
    $$
    \begin{aligned}
    \begin{bmatrix}
    v^{n, k+1} \\ 
    w^{n, k+1}
    \end{bmatrix} 
    =  
    \begin{bmatrix}
    \frac{1}{\Delta t} - 1.0 + \frac{(v^{n,k})^2}{3} & 1.0 \\
    -1.0 & \frac{\tau}{\Delta t} + b,
    \end{bmatrix}^{-1} 
    \begin{bmatrix}
    \frac{v^{n-1}}{\Delta t} + R I_{\text{ext}} \\
    a + \tau \frac{w^{n-1}}{\Delta t}
    \end{bmatrix}
    \end{aligned}.
    $$
    This approach will be used with tight error tolerances (\revNB{$10^{-10}$ relative increment}) to build up a ground truth against which to compare all the other schemes. A solution using $\Delta t = 10^{-3}$ is depicted in Figure \ref{fig:ground_truth_FNG}.

    \item Semi-implicit approach: a classical implicit explicit (IMEX) approach, neglecting any fixed-point iteration and decoupling the potential and gate equations:
    \begin{equation}
    \begin{aligned}
    \frac{v^{n} - v^{n-1}}{\Delta t} &= v^{n} - (v^{n-1})^2 v^{n}/3 - w^{n-1} + R I_{\text{ext}}, \\
    \tau\frac{w^{n} - w^{n-1}}{\Delta t} &= v^{n-1} + a - b w^{n-1}.
    \end{aligned}
    \end{equation}
    
    \item Accelerated implicitized approach: our proposal, using the acceleration scheme to enable fast computations of the explicit time-marching scheme within the Picard loop:
    \begin{equation}
    \begin{aligned}
    \frac{v^{n,k+1} - v^{n-1}}{\Delta t} &= v^{n,k} - (v^{n,k})^3/3 - w^{n,k} + R I_{\text{ext}}, \\
    \tau\frac{w^{n,k+1} - w^{n-1}}{\Delta t} &= v^{n,k} + a - b w^{n,k},
    \end{aligned}
    \end{equation}
    \revNB{where iterations are stopped when the same tolerance of the implicit approach is achieved.}

    \item Explicit approach: to compare with the previous approach, we check how the convergence speed improves due to AA compared to a purely explicit Picard subiterations.
\end{enumerate}

\begin{figure}
        \centering
        \includegraphics[width=0.9\linewidth]{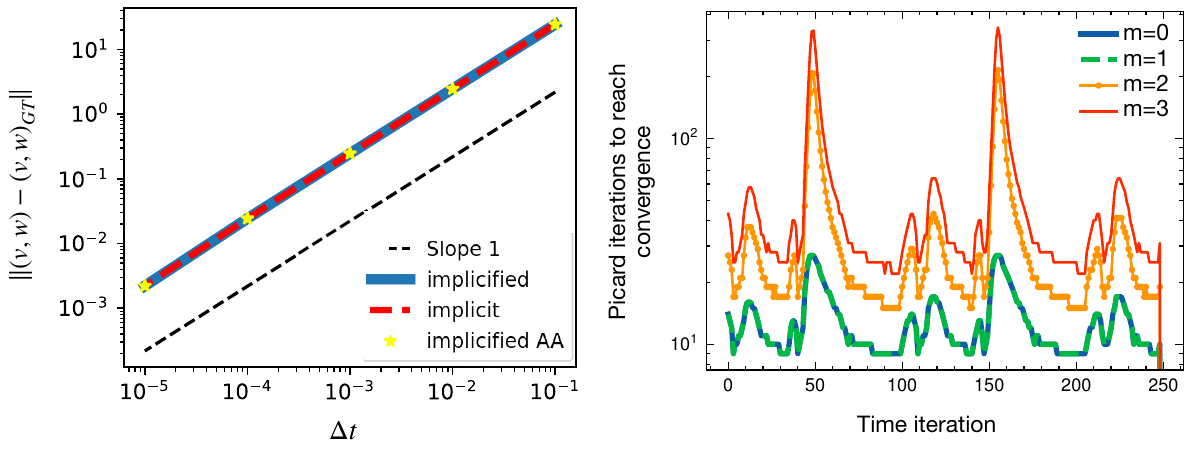}
        \caption{Convergence study (left) and sensitivity study to AA depth. }
        \label{fig:convergence_FN}
    \end{figure}

We can see a comparison of the numerical solutions using the aforementioned approaches in Figure \ref{fig:normalized_vw}. While the IMEX approach tends to largely offset the solution, both the Anderson accelerated approach (using $m=2$) and the implicitized approach, tend to perform better (and equal) to the implicit ground truth solution. Furthermore, we can see a convergence study in Figure \ref{fig:convergence_FN} (left), to verify how the implicitized approach (with the vanilla formulation and using AA) exhibits the same convergence rate as the fully implicit scheme.

    An additional study is depicted in Figure \ref{fig:convergence_FN} (right), where we test out the AA depth parameter $m$. In this admittedly simple numerical example, a depth larger than 2 leads to a minimal-norm Moore-Penrose Pseudoinverse, which we verify, in this particular case, to under-perform the vanilla implicitized approach.

\subsection{Application to heat equation}
In this section we show the results of using the formulation explicitly derived and analyzed in Section~\ref{sec:example-heat}. We start first by reporting the execution times of our implicitization iteration against a simple Picard solver. For this, we consider a diffusion coefficient $\mu=0.1$ and a varying range of timesteps to be able to test whether there is or not an improvement coming from acceleration. In addition we consider the following solvers to render the strategy HPC-enabled, all available in PETSc \citep{petsc-user-ref}:
    \begin{itemize}
        \item For the \ac{LS} problem in each acceleration step, we use the LSQR solver \citep{paige1982lsqr}.
        \item Each explicit solve of the mass matrix is done with the \ac{CG} algorithm and a Jacobi preconditioner \revNB{with a relative tolerance of $10^{-6}$, which can be setup with the PETSc options \texttt{-ksp\_type cg -pc\_type jacobi -ksp\_atol 0 -ksp\_rtol 1e-6.}}
    \end{itemize}
    We test this in a unit cube with 40 elements per side, which yields a total of roughly $70\,000$ \acp{dof} \revNB{, and iterate the fixed-point map given by the solution map of \eqref{eq:heat-imp-iter} (as in Section~\ref{sec:accel-imp}) until the relative norm of the k-th residual norm $|\vec R^k|^2_2 = \sum_l (R^k_l)^2$, given by $|\vec R^k| / |\vec R^0|$, is smaller than $10^{-6}$.}

    In Table~\ref{tab:heat-dofs} we report the average fixed point and \ac{CG} iterations during the first three instants. We note that, as expected from the convergence analysis, the timestep $\Delta t$ is subject to a CFL-like condition and cannot surpass a certain threshold. Acceleration increases the convergence radius as expected, and interestingly, we note that a larger window enables much larger timesteps. This goes against some commonly established evidence that suggests that larger depth is usually worse\revNB{, as theoretical convergence results might require the last $m$ iterates to be close to the fixed-point \citep{rebholz2023effect} and the least-squares system becomes asymptotically ill-posed due to collinearity of the residuals}. The \ac{CG} iterations are constant in practice, showing that our solution strategy for the mass matrix is adequate, and the increase in fixed-point iterations was expected from the increase in difficulty of the problem. \revNB{In the case in which the non-accelerated scheme converges, our approach reduces the number of fixed-point iterations by roughly a 96\%. We have additionally performed a numerical convergence test to verify that there is no deterioration in the solution when acceleration is considered. For this, we report in Table~\ref{tab:heat-convergence} the convergence rate obtained for a fixed timestep of $\Delta t=10^{-4}$ when doing uniform refinement for the explicit, implicit, and implicitized schemes for various acceleration depths. The errors are computed in the $L^\infty(H^1)$ norm, and we observe that the rates obtained are optimal and identical for all implicit solvers. This was to be expected because the solution computed is the same (up to the considered tolerance).} 

\begin{table}[ht!]
    \centering
    \begin{tabular}{r|r r r r}
         \toprule Depth ($m$) & $\Delta t=10^{-4}$ & $\Delta t=10^{-3}$ & $\Delta t=10^{-2}$ & $\Delta t=10^{-1}$ \\ \midrule
            0    & 180.7 (6.23) & --           & --     & --  \\
            10   & 7.7 (13.86)  & 28.0 (14.21) & --     & --  \\
            100  & 7.7 (13.86)  & 23.0 (14.61) & 73.0 (14.09) & 356.0 (14.98) \\ \bottomrule
    \end{tabular}
    \caption{Average fixed point iterations (\ac{CG} iterations) to solve the heat problem in the first 3 instants through an accelerated fixed-point iteration using varynig depth. Here '--' denotes non-convergence. }
    \label{tab:heat-dofs}
\end{table}

\revNB{
\begin{table}[h!]
\centering
\begin{tabular}{c|c|c|c|c|c|c}
\toprule
N & DoFs & Explicit & Implicit & $AA(2)$ & $AA(5)$ & $AA(10)$  \\ \midrule
4 & 125    &  --   & --   & --   & --   & --   \\
8 & 729    &  1.23 & 1.23 & 1.23 & 1.23 & 1.23 \\
16 & 4913  &  0.92 & 0.92 & 0.92 & 0.92 & 0.92 \\
32 & 35937 &  1.09 & 1.11 & 1.11 & 1.11 & 1.11 \\
\bottomrule
\end{tabular}
\caption{Convergence study for heat problem with $\Delta t \propto h$ for the Anderson Accelerated Picard scheme.}
\label{tab:heat-convergence}
\end{table}
}

We now show that this methodology can be easily modified to tackle nonlinear problems. As an example, we can modify our Picard iteration procedure to solve a heat equation with a p-Laplacian: 
    $$ \frac{u^{n, k+1}-u^{n-1}}{\Delta t} - 
    %\text{div}
    \nabla \cdot \left(|\nabla u^{n,k}|^{p-2}\nabla u^{n,k}\right) = 0 \qquad \text{in $(0,T)\times \Omega$}.$$
The problem for $p<2$ requires special treatment, see \citep{barnafi2025two} for a Levemberg-Marquardt formulation and \citep{loisel2020efficient} for a robust penalty method. We show the sensitivity of the iterations for increasing $p$ in Table~\ref{tab:heat-p}. We note that the method remains robust for increasing $p$, and present a mild and sublinear increase in iterations with the \acp{dof}.

\begin{table}[ht!]
    \centering
    \begin{tabular}{r|c c c c}
         \toprule \acp{dof} & $p=2$ & $p=3$ & $p=4$        & $p=5$  \\ \midrule
         216   & 1.67 & 1.67   & 1.67  & 1.67  \\
         1331  & 2.33 & 2.33   & 2.33  & 1.67  \\
         9261  & 3.67 & 3.33   & 3.33  & 3.00  \\
         68921 & 5.33 & 5.33   & 5.00  & 5.00    \\ \bottomrule
    \end{tabular}
    \caption{Average iterations to solve the nonlinear heat problem in the first 3 instants through an accelerated fixed-point iteration using a depth of $m=10$. }
    \label{tab:heat-p}
\end{table}

\subsection{Application to Navier-Stokes}

We conclude our numerical examples using the implicitized approach with the Navier-Stokes equations. So we aim to find $\vec u:\Omega \rightarrow \bR^2$ and $p: \Omega \rightarrow \bR$, such that:
\begin{equation}
\left\lbrace
\begin{aligned}
\rho \frac{\partial \bu}{\partial t} + \rho \bu \cdot \nabla \bu + \nabla p - \mu \Delta \bu  &= 0 \quad \text{in }  (0,T) \times \Omega ,\\
 \nabla \cdot \bu &= 0  \quad \text{in } (0,T) \times \Omega.
\end{aligned}
\right.
\label{eq:nse}
\end{equation}

\begin{figure}[!htbp]
    \centering
    \includegraphics[width=\linewidth]{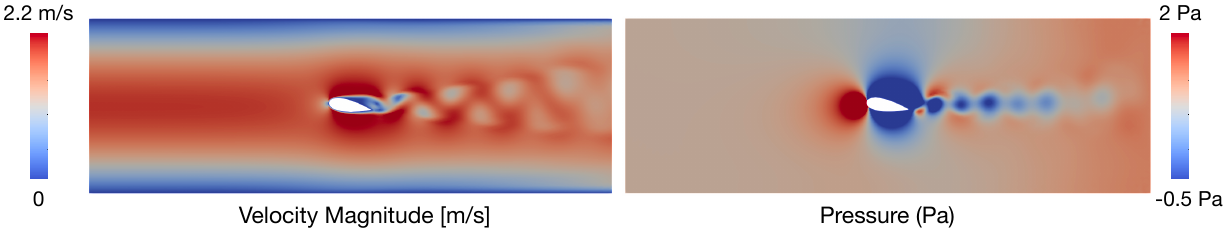}
    \caption{Numerical solutions of \eqref{eq:nse} using implicitized scheme with coarse $\Delta t = 0.15$ s and AA with $m=3$. We emphasize the success of the tandem implicitization-AA, as we observe how the solver diverges for $m=0$, numerically verifying the broader convergence radius when using AA \citep{lupo2019convergence}.}
    \label{fig:solutionNSE}
\end{figure}

Our working domain, depicted in Figure \ref{fig:solutionNSE}, is a two-dimensional symmetric NACA airfoil 0030. We repeat our tests to verify the feasibility of the implicitized scheme. We treat explicitly only the convective term, as this is better suited for large scale problems. In this way, we avoid matrix re-assembly within the Picard loop and still consider the incompressibility constraint. We set $\mu=0.001$ Pa s, $\rho = 1$ Kg/m$^3$, and a parabolic inlet velocity (with peak 2 m$/$s) that ensures the periodic vortex shedding in laminar regime. A traction-free condition is enforced in the domain right hand side, while no-slip boundary conditions are enforced at both, the NACA surface and the lateral walls. The flow starts with a zero-velocity initial condition.

The numerical solution of the governing dynamics is done by means of the finite element method, using the MAD library \citep{galarceThesis, MELLA2026108145, GALARCE202540}, built upon the linear algebra library PETSc \citep{petsc-user-ref}. We use an equal-order piece-wise linear formulation, adding SUPG-PSPG stabilization terms \citep{BROOKS1982199}. For our tests, we set the flexible GMRES as our linear solver, as it adapts well to variable Schur block preconditioning arising from inexact inner solvers within the lower and diagonal terms of the LDU factorization. We set a fixed-point error tolerance of $10^{-6}$ in the $H^1(\Omega) \times L^2(\Omega)$ natural norm, and the Picard loop is solved using FGMRES(100) (i.e. with a restart of 100) with an absolute $l^2$ error tolerance of $10^{-3}$. The block preconditioning is done with the built-in functions of PETSc, using loose inner solvers (with $l^2$ tolerances $5 \times 10^{-2}$ and $10^{-1}$ for the velocity and pressure blocks, respectively) to avoid computational bottlenecks due to unnecessary precision. We keep this strategy for all the numerical tests within this section. The computational mesh is adaptively refined in the NACA downstream region, leading to an unstructured triangular mesh with 17 000 nodes. Concerning time discretization, we choose to use an optimized second order backward difference formulation, as in \cite{OSSES2021114099, POBLETE2024125642}.

\begin{figure}[!htbp]
    \centering
    \includegraphics[width=0.9\linewidth]{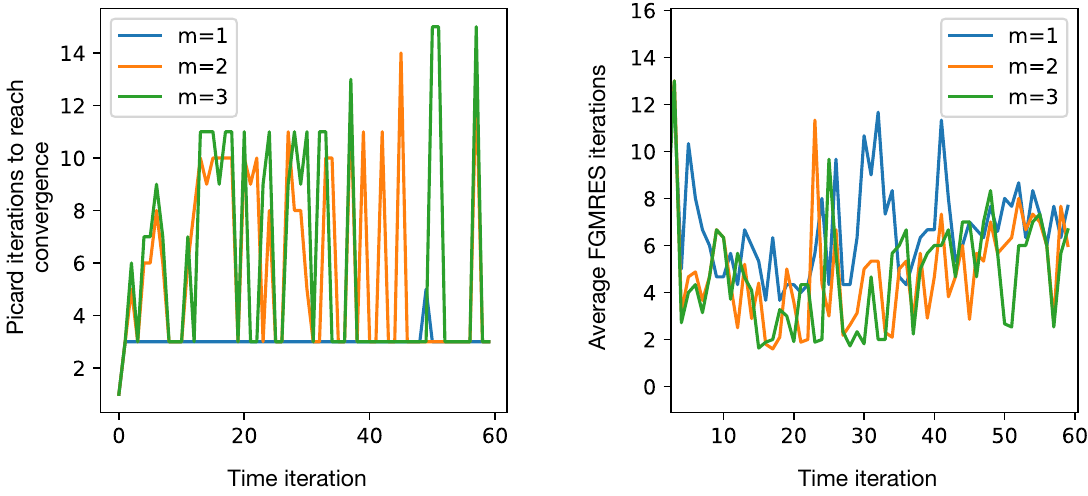}
    \caption{Picard iterations as time goes by. Optimal convergence speed is achieved in this numerical example for AA with $m=1$, and we highlight that the non-accelerated case $m=0$ fails to converge.}
    \label{fig:iterations_NSE}
\end{figure}

The first highlight in this numerical example relates to the fact that, with the coarse time-step we choose: $\Delta t = 0.15$ s, no convergence is reached if we turn off AA, while we verify how using any depth $m>0$ allow the convergence in a few Picard iterations. This is in agreement with results in the literature \citep{pollock2019anderson}. In addition, we observe that enlarging the AA depth brings a slower convergence (see Figure \ref{fig:iterations_NSE}, left) as has been reported for this equation \citep{pollock2019anderson}. We observe the best trade-off in terms of overshoot with respect to stability for $m=1$. In addition, the average amount of FGMRES iterations is depicted in Figure \ref{fig:iterations_NSE} (right), where we can observe that the slight decrease in average iterations for larger depths is related to faster solvers in the final Picard iterations. \revNB{We highlight that our tests have shown no advantage in using higher depth, and thus preferred not to report those solutions. This is in agreement with previous experience in solving Navier-Stokes \citep{pollock2019anderson}, and thus we show that this holds even when more relaxed Picard iterations are used for its solution.}

\begin{figure}[!htbp]
    \centering
    \includegraphics[width=0.4\linewidth]{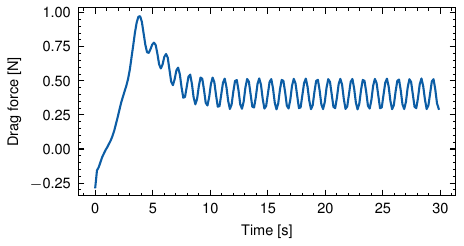}
    \caption{Drag force over the Naca profile obtained with $\Delta t = 0.15$ s, and $m=1$. The numerical solution exhibits stable outcomes for both the initial acceleration phase and the periodic fully developed regime.}
    \label{fig:drag_force}
\end{figure}

We additionally check the drag values obtained with this numerical approach. To compute the overall drag force over the obstacle, we use the variational approach described in \cite[Appendix D]{volker2016_femFlows} which is well suited for finite element solutions only using first order spatial derivatives. The resulting drag is depicted in Figure \ref{fig:drag_force}, showing a physically coherent behavior, similar to other benchmarks for thick airfoils found in the literature, in terms of vortex shedding frequency and flow separation \cite{nacaBen2021}.

\section{Conclusions and Outlook}

In this work we have elaborated on the simple idea that delayed terms in an explicit time discretization can be interpreted as a fixed-point iteration to solve the implicit nonlinear problem. This process is methodologically flexible as it does not make a difference regarding the linearity (or lack of it), and in many scenarios can perform similarly to state-of-the-art solvers. We believe this idea has been used implicitly in many contexts, such as for developing Picard solvers for multiphysics problems such as Navier-Stokes, Boussinesq, and the Magneto-hydro-dynamics equations, but this work establishes this methodology as a well-defined solver which is amenable to acceleration using \ac{AA}, which greatly improves the performance of the method. 

The aforementioned statements are supported by numerical evidence in a battery of examples, where we observe how our proposed pipeline reduces the number of fixed-point iterations up to $96 \%$ in diffusion problems, and how it enables the accurate resolution of non-linear fluid mechanic problems using very large time-steps ($\Delta t \sim 0.15$ s in a fully developed thick Naca airfoil Von-Kármán street), yielding a numerical scheme whose time step size is driven by the physical resolution that the user intends to capture, instead of cost-inefficient numerical constraints.

We believe that our approach is adequate to leverage different hardware such as GPU due to its low communication requirements at scale, and envision an exploration on other classical problems together with its extension to large scale simulations. 

%\section{Acknowledgments}

%N.A. Barnafi was supported by Centro de Modelamiento Matematico (CMM), Proyecto Basal FB210005, and by ANID Postdoctoral Proyecto 3230326. F. Galarce acknowledge the financial support of the Chilean National Agency for Research and Development (ANID) through the project Fondecyt No. 1250287. F. Galarce acknowledges DI Vinci PUCV Iniciación 039.731/2025 and the Horizon Europe - 2nd Opportunity OPPTY-MSCA/0125 funding from the Cyprus Research and Innovation Foundation. 

\bibliography{main}

\end{document}